\documentstyle[titlepage]{article}
\newtheorem{thm}{Theorem}[section]
\newtheorem{prop}[thm]{Proposition}
\newtheorem{cor}[thm]{Corollary}
\newtheorem{lem}[thm]{Lemma}
\newtheorem{conj}[thm]{Conjecture}
\newtheorem{exa}[thm]{Example}
\newtheorem{defn}[thm]{Definition}
\newcommand{\ben}{\begin{enumerate}}
\newcommand{\een}{\end{enumerate}}
\newcommand{\ble}{\begin{lem}}
\newcommand{\ele}{\end{lem}}
\newcommand{\bth}{\begin{thm}}
\newcommand{\eth}{\end{thm}}
\newcommand{\bpr}{\begin{prop}}
\newcommand{\epr}{\end{prop}}
\newcommand{\bco}{\begin{cor}}
\newcommand{\eco}{\end{cor}}
\newcommand{\bcon}{\begin{conj}}
\newcommand{\econ}{\end{conj}}
\newcommand{\bde}{\begin{defn}}
\newcommand{\ede}{\end{defn}}
\newcommand{\bex}{\begin{exa}}
\newcommand{\eex}{\end{exa}}
\newcommand{\barr}{\begin{array}}
\newcommand{\earr}{\end{array}}
\newcommand{\btab}{\begin{tabular}}
\newcommand{\etab}{\end{tabular}}
\newcommand{\beq}{\begin{equation}}
\newcommand{\eeq}{\end{equation}}
\newcommand{\bea}{\begin{eqnarray*}}
\newcommand{\eea}{\end{eqnarray*}}
\newcommand{\bce}{\begin{center}}
\newcommand{\ece}{\end{center}}
\newcommand{\bpi}{\begin{picture}}
\newcommand{\epi}{\end{picture}}
\newcommand{\bfi}{\begin{figure} \begin{center}}
\newcommand{\efi}{\end{center} \end{figure}}

\newcommand{\bsl}{\begin{slide}{}}
\newcommand{\esl}{\end{slide}}

\newcommand{\bib}{thebibliography}
\newcommand{\pf}{{\bf Proof.}}

\newcommand{\Qed}{\rule{1ex}{1ex} \medskip}

\newcommand{\emp}{\emptyset}

\newcommand{\fl}[1]{\lfloor #1 \rfloor}
\newcommand{\ce}[1]{\lceil #1 \rceil}

\newcommand{\la}{\lambda}

\newcommand{\bZ}{{\bf Z}}

\newcommand{\mod}{\mathop{\rm mod}}

\newcommand{\cho}{\choose}

\newcommand{\aim}{Adv. in Math.}

\newcommand{\jcta}{J. Combin. Theory Ser. A}

\newcommand{\tams}{Trans. Amer. Math. Soc.}

\begin{document}
\pagestyle{empty}
\title{Unimodality and the Reflection Principle}
\author{Bruce E. Sagan \\
Department of Mathematics \\ Michigan State University
\\ East Lansing, MI 48824-1027}

\date{\today \\[1in]
	\begin{flushleft}
	\end{flushleft}
       }
\maketitle

\begin{flushleft} Proposed running head: \end{flushleft}
	\begin{center} 
Unimodality and reflections
	\end{center}

Send proofs to:
\begin{center}
Bruce E. Sagan \\ Department of Mathematics \\Michigan State
University \\ East Lansing, MI 48824-1027\\[5pt]
Tel.: 517-355-8329\\
FAX: 517-336-1562\\
Email: sagan@mth.msu.edu
\end{center}

	\begin{abstract}
We show how lattice paths and the reflection principle can be used to
give easy proofs of unimodality results.  In particular, we give a
``one-line'' combinatorial proof of the unimodality of the binomial
coefficients.  Other examples include products of binomial
coefficients, polynomials related to the Legendre polynomials, and a
result connected to a conjecture of Simion.
	\end{abstract}
\pagestyle{plain}

\section{Introduction}

The Gessel-Viennot lattice path technique~\cite{gv:bdp,gv:dpp} has
proved useful in proving log concavity results~\cite{sag:lcs}.  The
purpose of this note is to show that lattice paths, and in particular
the Reflection Principle~\cite{moh:lpc}, are useful for demonstrating
the related property of unimodality.  One application is the simplest
combinatorial proof we know of that the binomial coefficients are
unimodal.  Such a proof is implicit in a standard use of the
Reflection Principle to show that the ballot numbers (which are
differences of binomial coefficients) are nonnegative, 
see~\cite[p. 95]{fel:ipt} or~\cite[p. 3]{moh:lpc}
However, the connection with unimodality does not seem to have been
made before.
This technique also applies to the four-step lattice paths
studied by DeTemple \& Robertson~\cite{dr:elf}, Cs\'aki, Mohanty \&
Saran~\cite{cms:rwp} , Guy, Krattenthaler \& Sagan~\cite{gks:lpr}, and
Beckenridge et. al.~\cite{bbcgnw:lpc}.  Lattice paths with diagonal
steps, which are related to Legendre polynomials~\cite{com:ac},
provide another example.  Finally, we prove a result similar to a
conjecture of Simion~\cite{sim:csn} and discuss some comments
and open questions.

Before beginning with the proofs, let us make some definitions.  A
sequence of real numbers $(a_k)_{k\geq0}$ is {\it unimodal} if there
is some index $m$ such that
	\beq						\label{uni}
	a_0\le a_1\le \ldots\le a_m\ge a_{m+1}\ge\ldots
	\eeq
Unimodal sequences arise in many areas of mathematics.  See Stanley's
survey article~\cite{sta:lcu} for details.  A related property is log
concavity.  A sequence is {\it log concave} if for all indices $i\ge1$ we
have
	$$a_{i-1}a_{i+1}\le a_i^2.$$
The following result is well known and easy to prove.
\bpr						\label{lcu}
Suppose $a_k>0$ for all $k$.  Then log concavity of $(a_k)_{k\geq0}$ 
implies unimodality.\hfill\Qed
\epr
In Section~\ref{fsp} we will see an example where unimodality
implies log concavity.

Let $\bZ^2$ denote the two-dimensional integer lattice.  A {\it
lattice path}, $p$, is a sequence $v_0,v_1,\ldots,v_n$ where $v_i\in\bZ^2$
for all $i$.  The integer $n$ is called the {\it length} of the path.
For example
	\beq							\label{p}
	p=(0,0),\ (0,1),\ (1,1),\ (2,1),\ (3,1),\ (3,2)
	\eeq
is a path of length five.
All of our paths will start at the origin, i.e., $v_0=(0,0)$.  Thus we
will often describe a path in terms of its steps, where the  {\it $i$th step}
is the vector $s_i$ from $v_{i-1}$ to $v_i$.  We will use square
brackets to enclose the coordinates of a vector so as to distinguish
it from an element of the lattice.  If $E=[1,0]$ and $N=[0,1]$ denote
the unit steps east and north, respectively, then the path~(\ref{p}) can
be written
	$$p=N,\ E,\ E,\ E,\ N.$$

Let $l$ be a line in the plane such that reflection in $l$ leaves
$\bZ^2$ invariant.  Consider a lattice path $p=v_0,\ldots,v_n$ that
intersects $l$ in at least one lattice point and let $v_k$ be the last
such point, i.e., the point with largest index.  Then the path
associated with $p$ via $l$ by {\it the reflection principle} is
	$$p'=v_0,\ldots,v_k,v'_{k+1},\ldots v'_n$$
where $v'_i$ is the reflection of $v_i$ in $l$.  By way of
illustration, if $p$ is the path in~(\ref{p}) and $l$ is the line
$y=x$ then
	$$p'=(0,0),\ (0,1),\ (1,1),\ (1,2),\ (1,3),\ (2,3).$$
We now have all the tools we will need to prove our unimodality
results.

\section{Two-step paths and binomial coefficients}	\label{tsp}

Let $T_v=T_{i,j}$ be the set of all lattice paths starting at the origin,
ending at $v=(i,j)\in\bZ^2$ and using only steps $N$ and $E$.  As is
well known, 
	$$|T_{i,j}|={i+j\choose i}$$
where $|\cdot|$ denotes cardinality.  So we can use such paths to
prove the unimodality of the binomial coefficients.  

\bth							\label{bin}
The sequence 
	$${n\choose0},{n\choose1},\ldots,{n\choose n}$$
is unimodal.
\eth
\pf\ By symmetry it suffices to find, for $k<n/2$, an injection
$T_v\hookrightarrow T_w$ where $v=(k,n-k),w=v+[1,-1]$.  Reflection
principle using 
the perpendicular bisector of the line segment $\overline{vw}$
provides such a map.\hfill\Qed

\section{Four-step paths}				\label{fsp}

Now consider lattice paths using steps in any one of the four directions
$N,E,S,W$ where $S=[0,-1]$ and $W=[-1,0]$.  Let $F_{i,j}(n)$ be the set
of all such paths that end at $(i,j)$ and have length $n$.  Note that
the length must be specified to make the set finite.  Furthermore, we
must have $n\equiv i+j(\mod 2)$ in order for this set to be nonempty,
and we will henceforth assume that this is the case.  Also let
$F^+_{i,j}(n)$ denote the subset of $F_{i,j}(n)$ consisting of all paths
that only use lattice points in the upper half-plane $y\geq0$.  The
cardinality of $F_{i,j}(n)$ was first derived in~\cite{dr:elf}.  A
combinatorial proof was given in~\cite{gks:lpr} where a formula for 
$|F^+_{i,j}(n)|$ was also obtained by the reflection principle.
\bpr							\label{FF+}
Let $r=(n-i-j)/2$ and $s=(n+i-j)/2$.  Then
	\bea
	|F_{i,j}(n)|&=&{n\choose r}{n\choose s},\\
	|F^+_{i,j}(n)|&=&{n\choose r}{n\choose s}-
		{n\choose r-1}{n\choose s-1}.\ \Qed
	\eea
\epr

The previous proposition can be used to prove unimodality of a number
of sequences involving binomial coefficients.  Two of the more
interesting follow.
\bth		
For any fixed integers $n,l$, the sequence
	\beq					\label{prod}
	{n\cho l}{n\cho0},{n\cho l-1}{n\choose1},\ldots,{n\cho0}{n\cho l}
	\eeq
is unimodal.
\eth  
\pf\ By Proposition~\ref{FF+}, if
$i=-l+2k$ and $j=n-l$ then 
	$$|F_{i,j}(n)|={n\cho l-k}{n\cho k}.$$
So by symmetry it suffices to find an 
an injection $F_{i,j}(n)\hookrightarrow F_{i+2,j}(n)$ when $i<0$.
Reflection principle using the line $y=i+1$ will do the
trick.\hfill\Qed

We note that this theorem could also be obtained by using the fact that
the binomial coefficient sequence is log concave, the 
result that a pointwise product of log concave sequences is log
concave, and Proposition~\ref{lcu}.  As another point of interest, it
is amusing 
to see that this unimodality result {\it implies} log concavity of the
binomial coefficients.  Indeed, when $l=2k$ then comparing terms in
the middle of~(\ref{prod}) yields
	$${n\cho k-1}{n\choose k+1}\le{n\cho k}^2.$$

In the next result it is convenient to interpret ${n\cho m}$ as zero
when $m<0$.
\bth
For any fixed integers $n,l$, the sequence
   $$\left\{{n\cho l-k}{n\cho k}-{n\cho l-k-1}{n\cho k-1}\right\}_{k\ge0}$$
is unimodal.
\eth  
\pf\  The proof is the same as the in the previous theorem, replacing
$F_{i,j}(n)$ by $F^+_{i,j}(n)$ throughout.\hfill\Qed

Although this can also be proved using a little bit of elementary
calculus, the proof is not as short or as elegant.

\section{Path with diagonal steps}			\label{pds}

We will now deal with paths having steps $N,E$ and $D$ where the last
is the diagonal step $D=[1,1]$.  Let $q$ be an indeterminate and
consider the polynomial
	$$D_{i,j}(q)=\sum_p q^{d(p)}$$
where the sum is over all such paths $p$ ending at $(i,j)$ and $d(p)$
is the number of diagonal steps on $p$.  For example,
$D_{2,4}(q)=6q^2+20q+15$.  Interest in $D_{i,j}(q)$ stems from the fact
that $D_{n,n}(1)=P_n(3)$ where $P_n(x)$ is a Legendre
polynomial~\cite{com:ac}.  These paths are
also related to a problem of I. Gessel~\cite{ges:pro} where he
essentially asked for a proof that $|D_{2n,2n+2}(-2)|$ is a Catalan
number.

Since we are now dealing with polynomials, we will need an associated
notion of unimodality.  Let $a(q)$ and $b(q)$ be two polynomials in
$q$.  Then we will write $a(q)\le_q b(q)$ if, for all $i$, the coefficient
of $q^i$ in $a(q)$ is less than or equal to the corresponding
coefficient in $b(q)$.  We now say that a sequence of polynomials,
$\{a_k(q)\}_{k\geq0}$ is $q$-unimodal if it satisfies~(\ref{uni})
with $\le$ replaced everywhere by $\le_q$ (and similarly for $\ge$).
For more information about $q$-unimodality, see~\cite{sag:lcs}.
\bth
The sequence
	\beq						\label{D}
	D_{0,n}(q),D_{1,n-1}(q),\ldots,D_{n,0}(q)
	\eeq
is $q$-unimodal.
\eth
\pf\  One can give the same proof as in Theorem~\ref{bin}.  All that
is needed is to note that reflection in the given line leaves the
number of $D$ steps on a given path invariant.\hfill\Qed

Actually this result can be seen as a corollary of the statement of
Theorem~\ref{bin}, rather than following from its proof,
because we have the following explicit formula for our polynomials
	$$D_{i,j}(q)=\sum_{d\ge0}{i+j-d\cho i-d,\ j-d,\ d}q^d.$$
To see this, note that a path to $(i,j)$ with $d$ diagonal steps must
have $i-d$ horizontal steps and $j-d$ vertical ones.  This explains
the trinomial coefficient in the sum.  Extracting the coefficient of
$q^d$ in the terms of the sequence~(\ref{D}), we obtain a sequence of
trinomial coefficients which have their third bottom index constant at
$d$.  Thus this sequence is unimodal by Theorem~\ref{bin}, and they
all reach their maximum at the same point, proving that~(\ref{D}) is
$q$-unimodal. 

\section{Simion's conjecture}				\label{sc}

A {\it partition} is a weakly decreasing sequence
$\la=(\la_0,\la_1,\ldots,\la_m)$ of positive integers.  It will be
convenient to also let $\la_i=0$ for $i>m$.

A path $p$ from $(0,0)$ to $(i,j)$ using only $N$ and $E$ steps must
stay inside the box $0\le x\le i,\ 0\le y\le j$.  We now consider paths
staying inside such a box with the Ferrers diagram of the partition
$\la$ removed from the upperleft corner.  To be precise, let
$T_{i,j}(\la)$ be the set of all $N,E$-paths $p$ to $(i,j)$ such that
any lattice point $(x,y)$ on $p$ satisfies $x\ge\la_{j-y}$.  For
example, if $\la=(1)$ then
	$$|T_{i,j}(\la)|={i+j\cho i}-1.$$
R. Simion~\cite[Conjecture 4.3]{sim:csn} has made the following conjecture.
\bcon
For any partition $\la$ and integer $n$, the sequence
	\beq
	\label{sim}
	|T_{0,n}(\la)|,\ |T_{1,n-1}(\la)|,\ \ldots,\ |T_{n,0}(\la)|
	\eeq
is unimodal.
\econ
Of course, when $\la=\emp$ this theorem reduces to Theorem~\ref{bin}.

Although we cannot prove this conjecture, we can get a related result.
Given $\la=(\la_0,\la_1,\ldots,\la_m)$ and $n$, define a sequence of
partitions by first 
removing rows from the Ferrers diagram of $\la$ and then adding
columns.  Specifically, for $0\le i\le\fl{n/2}$ let
	$$\la^i=(\la_i,\la_{i+1},\ldots,\la_m)$$
and
   $$\la^{\ce{n/2}+i}=(\la_{\ce{n/2}}+i,\la_{\ce{n/2}+1}+i,\ldots,\la_m+i)$$
where $\fl{\cdot}$ and $\ce{\cdot}$ are the round up and round down
functions, respectively.
\bth
For any partition $\la$ and integer $n$, the sequence
	$$|T_{0,n}(\la^0)|,\ |T_{1,n-1}(\la^1)|,\ \ldots,\ |T_{n,0}(\la^n)|$$
is unimodal.
\eth
\pf\  We will show that the first half of the sequence increases,
as the proof that the second half decreases is similar.  In fact, the
very same map used in the proof of Theorem~\ref{bin} will work.  One
need only note that if the path $p$ does not contain a point of
$\la^k$ then the image of $p$ will not contain a point of $\la^{k+1}$
since we reflect the portion of $p$ beyond the
{\it last} intersection with the line.\hfill\Qed

\section{Comments and open questions}

It would be interesting to find applications of this method to other
types of lattice paths.  One could also try using different lattices
and higher dimensional analogs.  One problem that needs to be overcome
in three or more dimensions is that reflection in a hyperplane which
is the perpendicular bisector of the line segment joining opposite
vertices of a hypercube does not stabilize the lattice of integer
points.  Hildebrand and Starkweather~\cite{hs:lci} have some results in this
direction.  

Although Simion's conjecture is still open, some special cases have
been resolved.  In particular, Hildebrand [personal communication] has
shown that the conjecture holds when $\la$ consists of a single row or
a single column.  He has also proved an asymptotic log concavity
result when $\la$ is a large rectangle.

\medskip

{\it Acknowledgement.}  I would like to thank the referee for helpful
suggestions.

\begin{\bib}{99}

\bibitem{bbcgnw:lpc} W. Breckenridge, P. Bos, G. Calvert, H.
Gastineau-Hills, A. Nelson and K. Wehrhahn, Lattice paths and
Catalan numbers, {\it Bulletin of the Institute of Combinatorics and
Applications} {\bf 1} (1991), 41--55.

\bibitem{com:ac} L. Comtet, ``Advanced Combinatorics,'' D. Reidel
Pub. Co., Dordrecht, 1974.

\bibitem{cms:rwp} E. Cs\'aki, S. G. Mohanty and S. Saran, On random
walks in a plane, {\it Ars Combin.} {\bf 29} (1990), 309--318.

\bibitem{dr:elf} D. W. DeTemple and J. M. Robertson, Equally likely
fixed length paths in graphs, {\it Ars Combin.} {\bf 17} (1984), 243--254.

\bibitem{fel:ipt} W. Feller, ``An Introduction to Probability Theory
and Its Applications, Vol. 1,'' Wiley, New York, NY, 1968.

\bibitem{ges:pro} I. Gessel, Problem number 10357, {\it Amer. Math.
Monthly} {\bf 101} (1994), 75.

\bibitem{gv:bdp} I.  Gessel  and  G.  Viennot,  Binomial 
determinants, paths, and hook length formulae, {\it \aim} {\bf 58}
(1985), 300--321.

\bibitem{gv:dpp} I. Gessel and G. Viennot, Determinants, paths, and
plane partitions, in preparation.

\bibitem{gks:lpr} R. K. Guy, C. Krattenthaler, and B. E. Sagan,
Lattice paths, reflections, and dimension-changing bijections, {\it
Ars Combin.} {\bf 34} (1992), 3--15.

\bibitem{hs:lci} M. Hildebrand and J. Starkweather, Log concavity
involving the number of paths from the origin to points along the line
$(a,b,c,d)+t(1,-1,1,-1)$, {\it Ars Combin.}, accepted.

\bibitem{moh:lpc} S. G. Mohanty, ``Lattice Path Counting
and Applications,'' Academic Press, New York, NY 1979.

\bibitem{sag:lcs} B. E. Sagan, Log concave sequences of symmetric
functions and analogs of the Jacobi-Trudi determinants, {\it \tams}
{\bf 329} (1992), 795--811.

\bibitem{sim:csn} R. Simion,  Combinatorial statistics on non-crossing
partitions, {\it \jcta} {\bf 66} (1994), 270--301.

\bibitem{sta:lcu} R. P. Stanley, Log-concave and unimodal 
sequences in algebra, combinatorics, and geometry, in
``Graph Theory and Its Applications: East and West,'' Ann. NY
Acad. Sci. {\bf 576} (1989), 500--535.

\end{\bib}

\end{document}